\documentclass[12pt,a4paper,notitlepage]{article}
\usepackage[english]{babel}
\usepackage{amsmath,amssymb,amsthm,amscd,latexsym}
\usepackage {graphicx}
\begin{document}
\title{Eighth lattice points}
\author{S. Adlaj\footnote{\textit{e-mail: SemjonAdlaj@gmail.com}}}
\date{May 30, 2011}
\parskip 6pt
\maketitle
\begin{abstract}
In a previous paper \cite{A}, a point of order 8 on an elliptic curve was calculated. Exploiting the well-known correspondence of the points on an elliptic curve with the points of a respective period parallelogram, we proceed to calculating all values at the eighth lattice points of an essential elliptic function which we, here, introduce.
\end{abstract}

\thispagestyle{empty}
\noindent
Given a complex parameter $\beta$, introduce an \textit{essential elliptic function} ${\cal R} = {\cal R}_\beta$ as the solution, with a (double) pole at zero, of the differential equation
$${\cal R}'^2 = 4 {\cal R} \left( {\cal R} + \beta \right) \left( {\cal R} + 1/\beta \right).$$
The function ${\cal R}$ differs by an additive constant from a Weierstrass elliptic function $\wp$ satisfying the differential equation
$$\wp'^2 = 4  \left( \wp - \alpha \right) \left( \wp - \alpha + \beta \right) \left( \wp - \alpha + 1 / \beta \right), \ \alpha:= \frac{\beta + 1/\beta}{3}.$$
Explicitly
$${\cal R} = \wp - \alpha.$$
Assuming that $\beta > 1$, the lattice $\Lambda$ of ${\cal R}$ (and of $\wp$) is rectangular \cite{WW}. Put
$$d = \beta - \frac{1}{\beta}, \ \delta = \sqrt{\frac{\beta}{d}}, \ \gamma^- = 1 - \frac{1}{\delta} < \gamma^+ = 1 + \frac{1}{\delta}, \ \gamma = 1 - \frac{i \beta}{\delta} \ \left( i:= \sqrt{-1} \ \right),$$
and set $\{ \beta_k \}$ and $\{ \delta_k \}$, $k = 0,1,2,3,$ to be the two sequences obtained by respectively replacing, with $\beta$ and $\delta$, the indeterminate $x$ in the ascending, for $x > 1$, sequence
$$x_0 = 1 - \sqrt{x + 1} < x_1 = 1 - \sqrt{1 + \frac{1}{x}} < x_2 = 1 + \sqrt{1 + \frac{1}{x}} < x_3 = 1 + \sqrt{x + 1}.$$
Put
$$\gamma_{02} = 1 + \frac{\delta_0 \delta_2}{\delta} < \gamma_{13} = 1 +\frac{\delta_1 \delta_3}{\delta} < \gamma_{01} = 1 + \frac{\delta_0 \delta_1}{\delta} < \gamma_{23} = 1 + \frac{\delta_2 \delta_3}{\delta},$$
$$\beta^- = 1 - \sqrt{1 - \frac{1}{\beta}} < \beta^+ = 1 + \sqrt{1 - \frac{1}{\beta}},$$
$$\delta^- = 1 - \frac{\sqrt{1 + \frac{1}{\beta}} - \sqrt{1 - \frac{1}{\beta}}}{\sqrt{2}} < \delta^+ = 1 + \frac{\sqrt{1 + \frac{1}{\beta}} - \sqrt{1 - \frac{1}{\beta}}}{\sqrt{2}},$$
$$\gamma_{0} = 1 - i \sqrt{\beta - 1}, \ \gamma_{-} = 1 - \frac{\delta^{-} \left( 1 + i \sqrt{\delta - 1} \right)}{\delta}, \ \gamma_{+} = 1 - \frac{\delta^{+} \left( 1 + i \sqrt{\delta - 1} \right)}{\delta}.$$
The values of ${\cal R}$ at the nodes of the lattice $\Lambda/8$ might now be represented as in the figure that follows (with the bar above a variable denoting complex conjugation). Yet, 4 values remain to be explicitely calculated. These are
$$\gamma_{41} = \left( 1 - \sqrt{\frac{\beta + 1}{2}} \ \right) \Bigg( 1 - \frac{1}{\beta} - i \sqrt{1 - \frac{1}{\beta^2}} + \frac{\sqrt{\beta + 1} \sqrt{\sqrt{\beta + 1} + \sqrt{2}}}{\beta}$$
$$\left( \sqrt{\frac{\beta}{\sqrt{\beta - 1}} - \frac{1}{\sqrt{\beta + 1}} + \sqrt{2}} - i \sqrt{\frac{\beta}{\sqrt{\beta - 1}} + \frac{1}{\sqrt{\beta + 1}} - \sqrt{2}} \ \right) \Bigg) - 1,$$
$$\gamma_{42} = \left( 1 - \sqrt{\frac{\beta + 1}{2}} \ \right) \Bigg( 1 - \frac{1}{\beta} + i \sqrt{1 - \frac{1}{\beta^2}} - \frac{\sqrt{\beta + 1} \sqrt{\sqrt{\beta + 1} + \sqrt{2}}}{\beta}$$
$$\left( \sqrt{\frac{\beta}{\sqrt{\beta - 1}} - \frac{1}{\sqrt{\beta + 1}} + \sqrt{2}} + i \sqrt{\frac{\beta}{\sqrt{\beta - 1}} + \frac{1}{\sqrt{\beta + 1}} - \sqrt{2}} \ \right) \Bigg) - 1,$$
$$\gamma_{43} = \left( \sqrt{\frac{\beta + 1}{2}} + 1 \right) \Bigg( 1 - \frac{1}{\beta} - i \sqrt{1 - \frac{1}{\beta^2}} - \frac{\sqrt{\beta + 1} \sqrt{\sqrt{\beta + 1} - \sqrt{2}}}{\beta}$$
$$\left( {\rm sgn}(\beta - b) \sqrt{\frac{\beta}{\sqrt{\beta - 1}} - \frac{1}{\sqrt{\beta + 1}}  - \sqrt{2}} - i \sqrt{\frac{\beta}{\sqrt{\beta - 1}} + \frac{1}{\sqrt{\beta + 1}} + \sqrt{2}} \ \right) \Bigg) - 1,$$
$$\gamma_{44} = \left( \sqrt{\frac{\beta + 1}{2}} + 1 \right) \Bigg( 1 - \frac{1}{\beta} + i \sqrt{1 - \frac{1}{\beta^2}} + \frac{\sqrt{\beta + 1} \sqrt{\sqrt{\beta + 1} - \sqrt{2}}}{\beta}$$
$$\left( {\rm sgn}(\beta - b) \sqrt{\frac{\beta}{\sqrt{\beta - 1}} - \frac{1}{\sqrt{\beta + 1}}  - \sqrt{2}} + i \sqrt{\frac{\beta}{\sqrt{\beta - 1}} + \frac{1}{\sqrt{\beta + 1}} + \sqrt{2}} \ \right) \Bigg) - 1,$$
where
$${\rm sgn}(x) :=
\begin{cases}
+1, & x > 0, \\
-1, & x < 0,
\end{cases}$$
and $b$ is the real root of the polynomial $x^3 - x^2 - x - 1$:
$$b = \frac{1 + \sqrt[3]{19 - 3 \sqrt{33}} +  \sqrt[3]{19 + 3 \sqrt{33}}}{3} \approx 1.839286759736968868.$$

\noindent
The mapping of the period rectangle by the function ${\cal R_\beta}$, for $\beta = \left( 3 + \sqrt{5} \right)/2$ $(\alpha = 1)$, is represented by the next (colorful) figure, where the red lines are images of lines parallel to the real axis, whereas the blue lines are images of lines parallel to the imaginary axis. In particular, the point at the center of one of the half-period rectangles is mapped to an intersection point of the red circle (centered at $-\beta$ with radius $\beta / \delta$) with the blue unit circle (centered at the origin). One might then observe that the value of the square of the logarithmic derivative of ${\cal R_\beta}$ at that point is $4 d$, coinciding (up to an independent of $\beta$ constant multiple) with the square root of the discriminant of the cubic on the right hand side of the differential equation satisfied by ${\cal R_\beta}$ (or by $\wp$). This value corresponds to the point of order 4 which one might obtain by doubling the point of order 8 described in \cite{A}. The points where ${\cal R_\beta}$ assumes the two values $\pm 1$ correspond to the points of order 4 given in \cite[\textsection 20.33, p. 444]{WW}. The latter observation provides an indirect, yet already convincing, a justification for introducing the essential elliptic function as an alternative, to Weierstrass elliptic function, for being more naturally, and thus more conveniently, investigated. The two circles are orthogonal and each circle inverts the lines of its own color from its inside to their corresponding same color lines outside and, of course, vice versa, whereas the lines whose color is that of the other circle (some of which continue beyond the boundaries of the figure) are inversion-invariant. Further elaboration might ensue in a subsequent paper.  

\begin{figure}[p]
\begin{picture}(560,500)(30,-64)
\put(58,-60){\parbox{404pt}{The values of the essential elliptic function ${\cal R}$, satisfying the differential equation}}
\put(112,-80){\parbox{400pt}{${\cal R}'^2 = 4 {\cal R} \left( {\cal R} + \beta \right) \left( {\cal R} + 1 / \beta \right)$, for $\beta > 1$, whose lattice is $\Lambda$,}}
\put(172,-100){\parbox{400pt}{at the nodes of the lattice $\Lambda$/8.}}
\fontsize{10}{0}
\put(18,20){\circle{28}}
\put(18,20){\circle{32}}
\put(258,20){\circle{32}}
\put(498,20){\circle{28}}
\put(498,20){\circle{32}}
\put(18,260){\circle{32}}
\put(258,260){\circle{32}}
\put(498,260){\circle{32}}
\put(18,500){\circle{28}}
\put(18,500){\circle{32}}
\put(258,500){\circle{32}}
\put(498,500){\circle{28}}
\put(498,500){\circle{32}}
\put(18,20){\circle{36}}
\put(258,20){\circle{36}}
\put(498,20){\circle{36}}
\put(18,260){\circle{36}}
\put(258,260){\circle{36}}
\put(498,260){\circle{36}}
\put(18,500){\circle{36}}
\put(258,500){\circle{36}}
\put(498,500){\circle{36}}

\put(138,20){\circle{36}}
\put(378,20){\circle{36}}
\put(18,140){\circle{36}}
\put(138,140){\circle{36}}
\put(258,140){\circle{36}}
\put(378,140){\circle{36}}
\put(498,140){\circle{36}}
\put(378,260){\circle{36}}
\put(138,260){\circle{36}}
\put(18,380){\circle{36}}
\put(498,380){\circle{36}}
\put(138,380){\circle{36}}
\put(258,380){\circle{36}}
\put(378,380){\circle{36}}
\put(138,500){\circle{36}}
\put(378,500){\circle{36}}

\put(18,20){\circle{40}}
\put(78,20){\circle{40}}
\put(138,20){\circle{40}}
\put(198,20){\circle{40}}
\put(258,20){\circle{40}}
\put(318,20){\circle{40}}   
\put(378,20){\circle{40}}
\put(438,20){\circle{40}}
\put(498,20){\circle{40}}
\put(18,80){\circle{40}}
\put(78,80){\circle{40}}
\put(138,80){\circle{40}}
\put(198,80){\circle{40}}
\put(258,80){\circle{40}}
\put(318,80){\circle{40}}   
\put(378,80){\circle{40}}
\put(438,80){\circle{40}}
\put(498,80){\circle{40}}
\put(18,140){\circle{40}}
\put(78,140){\circle{40}}
\put(138,140){\circle{40}}
\put(198,140){\circle{40}}
\put(258,140){\circle{40}}
\put(318,140){\circle{40}}   
\put(378,140){\circle{40}}
\put(438,140){\circle{40}}
\put(498,140){\circle{40}}
\put(18,200){\circle{40}}
\put(78,200){\circle{40}}
\put(138,200){\circle{40}}
\put(198,200){\circle{40}}
\put(258,200){\circle{40}}
\put(318,200){\circle{40}}   
\put(378,200){\circle{40}}
\put(438,200){\circle{40}}
\put(498,200){\circle{40}}
\put(18,260){\circle{40}}
\put(78,260){\circle{40}}
\put(138,260){\circle{40}}
\put(198,260){\circle{40}}
\put(258,260){\circle{40}}
\put(318,260){\circle{40}}   
\put(378,260){\circle{40}}
\put(438,260){\circle{40}}
\put(498,260){\circle{40}}
\put(18,320){\circle{40}}
\put(78,320){\circle{40}}
\put(138,320){\circle{40}}
\put(198,320){\circle{40}}
\put(258,320){\circle{40}}
\put(318,320){\circle{40}}   
\put(378,320){\circle{40}}
\put(438,320){\circle{40}}
\put(498,320){\circle{40}}
\put(18,380){\circle{40}}
\put(78,380){\circle{40}}
\put(138,380){\circle{40}}
\put(198,380){\circle{40}}
\put(258,380){\circle{40}}
\put(318,380){\circle{40}}   
\put(378,380){\circle{40}}
\put(438,380){\circle{40}}
\put(498,380){\circle{40}}
\put(18,440){\circle{40}}
\put(78,440){\circle{40}}
\put(138,440){\circle{40}}
\put(198,440){\circle{40}}
\put(258,440){\circle{40}}
\put(318,440){\circle{40}}   
\put(378,440){\circle{40}}
\put(438,440){\circle{40}}
\put(498,440){\circle{40}}
\put(18,500){\circle{40}}
\put(78,500){\circle{40}}
\put(138,500){\circle{40}}
\put(198,500){\circle{40}}
\put(258,500){\circle{40}}
\put(318,500){\circle{40}}   
\put(378,500){\circle{40}}
\put(438,500){\circle{40}}
\put(498,500){\circle{40}}

\put(38,20){\line(1,0){20}}
\put(98,20){\line(1,0){20}}
\put(158,20){\line(1,0){20}}
\put(218,20){\line(1,0){20}}
\put(278,20){\line(1,0){20}}
\put(338,20){\line(1,0){20}}
\put(398,20){\line(1,0){20}}
\put(458,20){\line(1,0){20}}
\put(18,40){\line(0,1){20}}
\put(78,40){\line(0,1){20}}
\put(138,40){\line(0,1){20}}
\put(198,40){\line(0,1){20}}
\put(258,40){\line(0,1){20}}
\put(318,40){\line(0,1){20}}
\put(378,40){\line(0,1){20}}
\put(438,40){\line(0,1){20}}
\put(498,40){\line(0,1){20}}

\put(38,80){\line(1,0){20}}
\put(98,80){\line(1,0){20}}
\put(158,80){\line(1,0){20}}
\put(218,80){\line(1,0){20}}
\put(278,80){\line(1,0){20}}
\put(338,80){\line(1,0){20}}
\put(398,80){\line(1,0){20}}
\put(458,80){\line(1,0){20}}
\put(18,100){\line(0,1){20}}
\put(78,100){\line(0,1){20}}
\put(138,100){\line(0,1){20}}
\put(198,100){\line(0,1){20}}
\put(258,100){\line(0,1){20}}
\put(318,100){\line(0,1){20}}
\put(378,100){\line(0,1){20}}
\put(438,100){\line(0,1){20}}
\put(498,100){\line(0,1){20}}

\put(38,140){\line(1,0){20}}
\put(98,140){\line(1,0){20}}
\put(158,140){\line(1,0){20}}
\put(218,140){\line(1,0){20}}
\put(278,140){\line(1,0){20}}
\put(338,140){\line(1,0){20}}
\put(398,140){\line(1,0){20}}
\put(458,140){\line(1,0){20}}
\put(18,160){\line(0,1){20}}
\put(78,160){\line(0,1){20}}
\put(138,160){\line(0,1){20}}
\put(198,160){\line(0,1){20}}
\put(258,160){\line(0,1){20}}
\put(318,160){\line(0,1){20}}
\put(378,160){\line(0,1){20}}
\put(438,160){\line(0,1){20}}
\put(498,160){\line(0,1){20}}

\put(38,200){\line(1,0){20}}
\put(98,200){\line(1,0){20}}
\put(158,200){\line(1,0){20}}
\put(218,200){\line(1,0){20}}
\put(278,200){\line(1,0){20}}
\put(338,200){\line(1,0){20}}
\put(398,200){\line(1,0){20}}
\put(458,200){\line(1,0){20}}
\put(18,220){\line(0,1){20}}
\put(78,220){\line(0,1){20}}
\put(138,220){\line(0,1){20}}
\put(198,220){\line(0,1){20}}
\put(258,220){\line(0,1){20}}
\put(318,220){\line(0,1){20}}
\put(378,220){\line(0,1){20}}
\put(438,220){\line(0,1){20}}
\put(498,220){\line(0,1){20}}

\put(38,260){\line(1,0){20}}
\put(98,260){\line(1,0){20}}
\put(158,260){\line(1,0){20}}
\put(218,260){\line(1,0){20}}
\put(278,260){\line(1,0){20}}
\put(338,260){\line(1,0){20}}
\put(398,260){\line(1,0){20}}
\put(458,260){\line(1,0){20}}
\put(18,280){\line(0,1){20}}
\put(78,280){\line(0,1){20}}
\put(138,280){\line(0,1){20}}
\put(198,280){\line(0,1){20}}
\put(258,280){\line(0,1){20}}
\put(318,280){\line(0,1){20}}
\put(378,280){\line(0,1){20}}
\put(438,280){\line(0,1){20}}
\put(498,280){\line(0,1){20}}

\put(38,320){\line(1,0){20}}
\put(98,320){\line(1,0){20}}
\put(158,320){\line(1,0){20}}
\put(218,320){\line(1,0){20}}
\put(278,320){\line(1,0){20}}
\put(338,320){\line(1,0){20}}
\put(398,320){\line(1,0){20}}
\put(458,320){\line(1,0){20}}
\put(18,340){\line(0,1){20}}
\put(78,340){\line(0,1){20}}
\put(138,340){\line(0,1){20}}
\put(198,340){\line(0,1){20}}
\put(258,340){\line(0,1){20}}
\put(318,340){\line(0,1){20}}
\put(378,340){\line(0,1){20}}
\put(438,340){\line(0,1){20}}
\put(498,340){\line(0,1){20}}

\put(38,380){\line(1,0){20}}
\put(98,380){\line(1,0){20}}
\put(158,380){\line(1,0){20}}
\put(218,380){\line(1,0){20}}
\put(278,380){\line(1,0){20}}
\put(338,380){\line(1,0){20}}
\put(398,380){\line(1,0){20}}
\put(458,380){\line(1,0){20}}
\put(18,400){\line(0,1){20}}
\put(78,400){\line(0,1){20}}
\put(138,400){\line(0,1){20}}
\put(198,400){\line(0,1){20}}
\put(258,400){\line(0,1){20}}
\put(318,400){\line(0,1){20}}
\put(378,400){\line(0,1){20}}
\put(438,400){\line(0,1){20}}
\put(498,400){\line(0,1){20}}

\put(38,440){\line(1,0){20}}
\put(98,440){\line(1,0){20}}
\put(158,440){\line(1,0){20}}
\put(218,440){\line(1,0){20}}
\put(278,440){\line(1,0){20}}
\put(338,440){\line(1,0){20}}
\put(398,440){\line(1,0){20}}
\put(458,440){\line(1,0){20}}
\put(18,460){\line(0,1){20}}
\put(78,460){\line(0,1){20}}
\put(138,460){\line(0,1){20}}
\put(198,460){\line(0,1){20}}
\put(258,460){\line(0,1){20}}
\put(318,460){\line(0,1){20}}
\put(378,460){\line(0,1){20}}
\put(438,460){\line(0,1){20}}
\put(498,460){\line(0,1){20}}

\put(38,500){\line(1,0){20}}
\put(98,500){\line(1,0){20}}
\put(158,500){\line(1,0){20}}
\put(218,500){\line(1,0){20}}
\put(278,500){\line(1,0){20}}
\put(338,500){\line(1,0){20}}
\put(398,500){\line(1,0){20}}
\put(458,500){\line(1,0){20}}

\put(12,16){$\infty$}
\put(136,16){$1$}
\put(255,16){$0$}
\put(376,16){$1$}
\put(493,16){$\infty$}

\put(5,136){$- \beta \hspace{.02cm} \gamma^+$}
\put(126,136){$\displaystyle - \frac{\overline{\gamma}}{\beta}$}
\put(243,136){$\displaystyle - \beta \hspace{.02cm} \gamma^-$}
\put(366,136){$\displaystyle - \frac{\gamma}{\beta}$}
\put(484,136){$- \beta \hspace{.02cm} \gamma^+$}

\put(8,256){$- \beta$}
\put(128,256){$- 1$}
\put(247,256){$\displaystyle - \frac{1}{\beta}$}
\put(368,256){$- 1$}
\put(488,256){$- \beta$}

\put(5,376){$- \beta \hspace{.02cm} \gamma^+$}
\put(126,376){$\displaystyle - \frac{\gamma}{\beta}$}
\put(243,376){$\displaystyle - \beta \hspace{.02cm} \gamma^-$}
\put(366,376){$\displaystyle - \frac{\overline{\gamma}}{\beta}$}
\put(484,376){$- \beta \hspace{.02cm} \gamma^+$}

\put(12,496){$\infty$}
\put(134,496){$1$}
\put(255,496){$0$}
\put(374,496){$1$}
\put(493,496){$\infty$}

\put(68,16){$\beta_2 \beta_3$}
\put(428,16){$\beta_2 \beta_3$}
\put(68,496){$\beta_2 \beta_3$}
\put(428,496){$\beta_2 \beta_3$}
\put(188,16){$\beta_0 \beta_1$}
\put(308,16){$\beta_0 \beta_1$}
\put(188,496){$\beta_0 \beta_1$}
\put(308,496){$\beta_0 \beta_1$}

\put(68,256){$\beta_0 \beta_2$}
\put(188,256){$\beta_1 \beta_3$}
\put(308,256){$\beta_1 \beta_3$}
\put(428,256){$\beta_0 \beta_2$}

\put(4,436){$- \beta \gamma_{23}$}
\put(124,436){$- \beta \hspace{.02cm} \gamma_{+}$}
\put(242,436){$- \beta \gamma_{02}$}
\put(364,436){$- \beta \hspace{.04cm} \overline{\gamma}_{+}$}
\put(484,436){$- \beta \gamma_{23}$}

\put(4,316){$- \beta \gamma_{01}$}
\put(124,316){$- \beta \hspace{.02cm} \gamma_{-}$}
\put(242,316){$- \beta \gamma_{13}$}
\put(362,316){$- \beta \hspace{.04cm} \overline{\gamma}_{-}$}
\put(484,316){$- \beta \gamma_{01}$}

\put(4,196){$- \beta \gamma_{01}$}
\put(124,196){$- \beta \hspace{.04cm} \overline{\gamma}_{-}$}
\put(242,196){$- \beta \gamma_{13}$}
\put(364,196){$- \beta \hspace{.02cm} \gamma_{-}$}
\put(484,196){$- \beta \gamma_{01}$}

\put(4,76){$- \beta \gamma_{23}$}
\put(124,76){$- \beta \hspace{.04cm} \overline{\gamma}_{+}$}
\put(242,76){$- \beta \gamma_{02}$}
\put(364,76){$- \beta \hspace{.02cm} \gamma_{+}$}
\put(484,76){$- \beta \gamma_{23}$}

\put(62,376){$- \beta^+ \gamma_{0}$}
\put(182,376){$- \beta^- \gamma_{0}$}
\put(302,376){$- \beta^- \overline{\gamma}_{0}$}
\put(422,376){$- \beta^+ \overline{\gamma}_{0}$}
\put(62,136){$- \beta^+ \overline{\gamma}_{0}$}
\put(182,136){$- \beta^- \overline{\gamma}_{0}$}
\put(302,136){$- \beta^- \gamma_{0}$}
\put(422,136){$- \beta^+ \gamma_{0}$}

\put(72,436){$\gamma_{44}$}
\put(192,436){$\gamma_{43}$}
\put(312,316){$\overline{\gamma}_{42}$}
\put(432,316){$\overline{\gamma}_{41}$}
\put(72,196){$\overline{\gamma}_{41}$}
\put(192,196){$\overline{\gamma}_{42}$}
\put(312,76){$\gamma_{43}$}
\put(432,76){$\gamma_{44}$}
\put(72,316){$\gamma_{41}$}
\put(192,316){$\gamma_{42}$}
\put(312,436){$\overline{\gamma}_{43}$}
\put(432,436){$\overline{\gamma}_{44}$}
\put(72,76){$\overline{\gamma}_{44}$}
\put(192,76){$\overline{\gamma}_{43}$}
\put(312,196){$\gamma_{42}$}
\put(432,196){$\gamma_{41}$}
\end{picture}
\end{figure}

\pagestyle{empty}
\begin{figure}[p]
\begin{picture}(560,720)(40,-40)
\put(-17,372){\parbox{400pt}{$- \beta \gamma^+$}}
\put(151,372){\parbox{400pt}{$- \beta \gamma_{01}$}}
\put(186,372){\parbox{400pt}{$- \beta$}}
\put(360,372){\parbox{400pt}{$\beta_1 \beta_3$}}
\put(248,372){\parbox{400pt}{$\beta_0 \beta_2$}}
\put(331,373){\parbox{400pt}{$-1$}}
\put(427,374){\parbox{400pt}{$0$}}
\put(438,372){\parbox{400pt}{$\beta_0 \beta_1$}}
\put(510,374){\parbox{400pt}{$1$}}
\put(378,469){\parbox{400pt}{$- \frac{\gamma}{\beta}$}}
\put(339,429){\parbox{400pt}{$- \beta \gamma_-$}}
\put(459,462){\parbox{400pt}{$- \beta \gamma_+$}}
\put(256,583){\parbox{400pt}{$- \beta^+ \gamma_0$}}
\put(388,410){\parbox{400pt}{$- \beta^- \gamma_0$}}
\put(240,456){\parbox{400pt}{$\gamma_{41}$}}
\put(431,398){\parbox{400pt}{$\gamma_{43}$}}
\put(490,772){\parbox{400pt}{$\gamma_{44}$}}
\includegraphics[width=180mm]{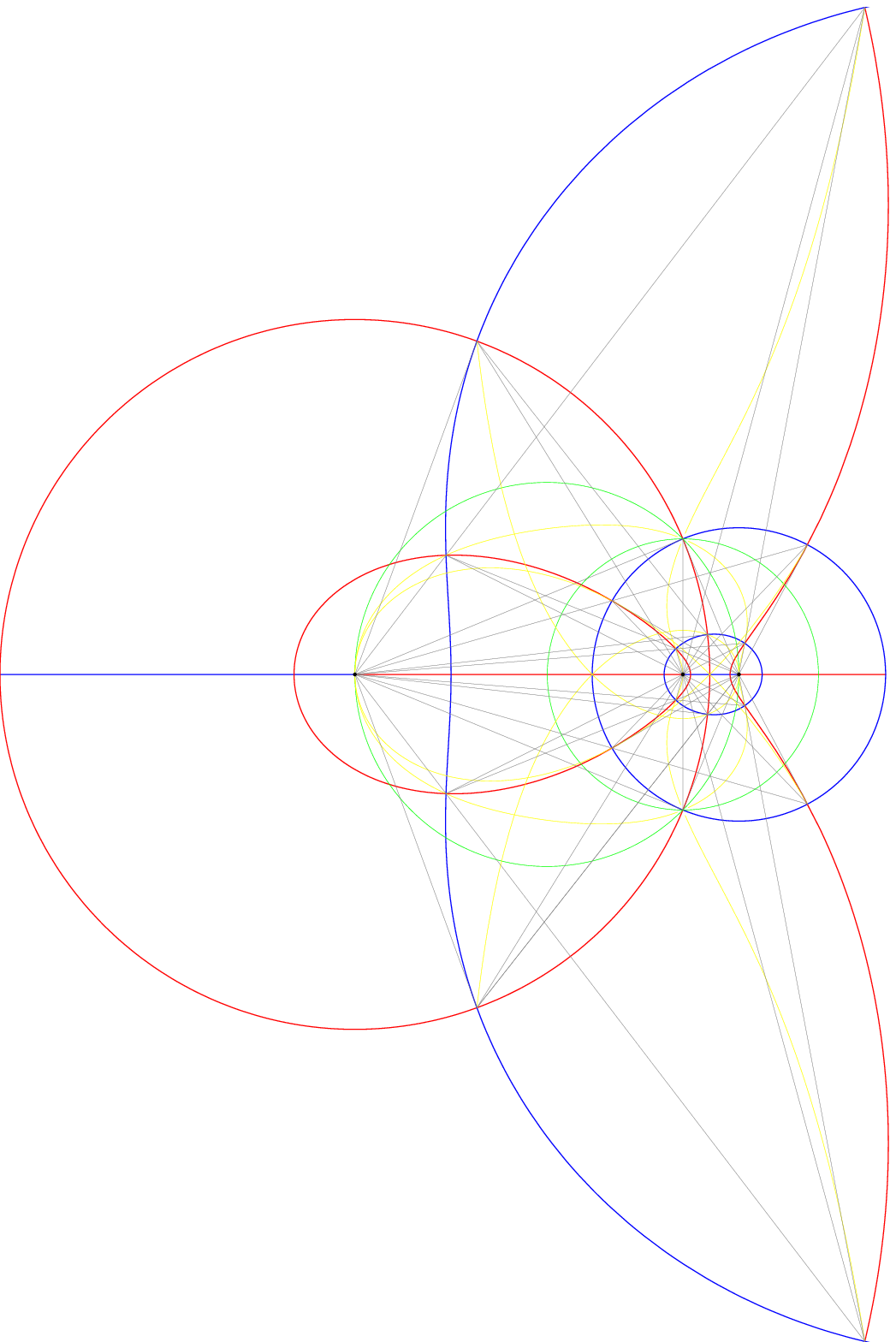}
\end{picture}
\end{figure}

\end{document}